\documentclass[a4paper,12pt]{article}

\usepackage[english]{babel}
\usepackage[letterpaper,top=2cm,bottom=2cm,left=3cm,right=3cm,marginparwidth=1.75cm]{geometry}
\usepackage[utf8]{inputenc}
\usepackage[T1]{fontenc}
\usepackage{url}
\usepackage{amsmath,amssymb,amsthm,amsfonts}
\usepackage{graphicx}
\usepackage[colorlinks=true, allcolors=red]{hyperref}

\theoremstyle{plain}
\newtheorem{theorem}{Theorem}[section]
\newtheorem{lemma}[theorem]{Lemma}

\newcommand{\Q}{\mathbb{Q}}
\newcommand{\Z}{\mathbb{Z}}
\newcommand{\bI}{\mathbb{L}}
\newcommand{\ol}[1]{\overline{#1}}
\numberwithin{equation}{section}
\allowdisplaybreaks
\begin{document}

\begin{center}
{\bf {\Large Padovan and Perrin Numbers as\\[4pt] Sums of Two Lucas Numbers}}
\vspace{5mm}

{\large O. Salah$^{*1}$, A. Elsonbaty$^{*2}$, M. Anwar$^{*3}$ and M. A. Seoud$^{*4}$}
\vspace{3mm}

$^*$Department of Mathematics, Faculty of Science, Ain Shams University\\
e-mail: osamasalah@sci.asu.edu.eg$^{*1}$,\\
ahmedelsonbaty@sci.asu.edu.eg$^{*2}$,\\
mohamedanwar@sci.asu.edu.eg$^{*3}$,\\
m.a.seoud@sci.asu.edu.eg$^{*4}$.
\end{center}

\begin{abstract}
Let $\{P_\ell\}_{\ell\ge0}$, $\{R_\ell\}_{\ell\ge0}$ and $\{L_n\}_{n\ge0}$ be the sequences of Padovan, Perrin and Lucas numbers, respectively. In this paper we determine all Padovan numbers and all Perrin numbers that can be written as a sum of two Lucas numbers, that is, all non-negative solutions $(\ell,r,s)$ with $r\ge s$ of the Diophantine equations $P_\ell=L_r+L_s$ and $R_\ell=L_r+L_s$. Our proofs combine Baker's method for linear forms in logarithms (in Matveev's version), the Dujella--Peth\H{o} reduction method, and a continued fraction argument based on Legendre's theorem.
\end{abstract}
\begin{flushleft}
\end{flushleft}
\begin{flushleft}
\end{flushleft}

\section{Introduction}

The Padovan sequence $\{P_\ell\}_{\ell\ge0}$ is defined by $P_\ell=P_{\ell-2}+P_{\ell-3}$ for $\ell\ge3$ with $P_0=P_1=P_2=1$. R. Padovan introduced it as a sequence connected to the plastic number previously defined by H. V. D. Laan. The sequence of Perrin numbers $\{R_\ell\}_{\ell\ge0}$ obeys the same recurrence $R_\ell=R_{\ell-2}+R_{\ell-3}$ for $\ell\ge3$, but begins with $R_0=3$, $R_1=0$, $R_2=2$. The Padovan and Perrin numbers therefore share the characteristic polynomial and many properties. Diophantine equations involving these sequences have been studied by many mathematicians; see for instance \cite{b1,b2,b3,b4,b6,b7,b8,b10,b11,b12}. The sequence of Lucas numbers $\{L_n\}_{n\ge0}$ is defined by $L_0=2$, $L_1=1$ and $L_n=L_{n-1}+L_{n-2}$ for $n\ge2$.

Closely related problems were solved recently. All Padovan and Perrin numbers which are sums of two Jacobsthal numbers were found in \cite{b15}; the authors treated Lucas versus Jacobsthal numbers in \cite{b13}; and the same tools were used for Narayana numbers in \cite{b14}. Here we replace the Jacobsthal sequence of \cite{b15} by the Lucas sequence, and we study the two Diophantine equations
\begin{equation}\label{eq:padmain}
P_\ell=L_r+L_s,\qquad r\ge s,
\end{equation}
and
\begin{equation}\label{eq:permain}
R_\ell=L_r+L_s,\qquad r\ge s,
\end{equation}
that is, we look for all Padovan numbers and all Perrin numbers which are sums of two Lucas numbers. We use linear forms in logarithms and reduction procedures. The complete sets of solutions are given by the following two theorems.

\begin{theorem}\label{thm:pad}
The only non-negative solutions $(\ell,r,s)$ which satisfy Eq.\eqref{eq:padmain} with $r\ge s$ belong to the set
\[
\left\{
\begin{array}{l}
(3,1,1),\ (4,1,1),\ (5,1,0),\ (6,0,0),\ (6,2,1),\\
(7,2,0),\ (7,3,1),\ (8,3,2),\ (9,4,0),\ (10,5,1),\\
(12,6,2),\ (15,8,0),\ (16,8,6),\ (20,11,1),\ (22,12,7)
\end{array}
\right\}.
\]
\end{theorem}

\begin{theorem}\label{thm:per}
The only non-negative solutions $(\ell,r,s)$ which satisfy Eq.\eqref{eq:permain} with $r\ge s$ belong to the set
\[
\left\{
\begin{array}{llllll}
(0,1,0), & (2,1,1), & (3,1,0), & (4,1,1), & (5,2,0), & (5,3,1),\\
(6,2,0), & (6,3,1), & (7,3,2), & (8,4,2), & (9,5,1), & (11,5,5),\\
(11,6,3), & (12,6,5), & (14,8,3), & (23,12,12), & (23,13,10). &
\end{array}
\right\}
\]
\end{theorem}

\emph{Organisation.} Section~2 recalls the facts we need on the three sequences, on heights, and the tools of Matveev, Dujella--Peth\H{o} and Legendre. Section~3 gives an estimate for $\ell$ in terms of $r$ valid for both equations, and Sections~4 and~5 prove Theorems~\ref{thm:pad} and~\ref{thm:per}, respectively.

\section{Preliminary results}

\subsection{Lucas sequence}
Let $\theta=\dfrac{1+\sqrt5}{2}$ be the golden ratio and $\eta=\dfrac{1-\sqrt5}{2}=-\theta^{-1}$ its conjugate, so that $|\eta|<1$. The Binet formula for the Lucas sequence is
\begin{equation}\label{eq:binetL}
L_n=\theta^{n}+\eta^{n}\qquad\text{for all }n\ge0 .
\end{equation}
By induction we have
\begin{equation}\label{eq:Lbd}
\theta^{n-1}\le L_n\le 2\theta^{n}\qquad\text{for all }n\ge0 .
\end{equation}
The cases $n=0,1$ are immediate, and the step from $n,n+1$ to $n+2$ follows from $L_{n+2}=L_{n+1}+L_n$ and $\theta^{2}=\theta+1$.

\subsection{Padovan sequence}
The characteristic equation of the Padovan and of the Perrin sequence is
\[
\Psi(x):=x^{3}-x-1=0,
\]
which has one real root $\alpha$ and two complex conjugate roots $\beta$ and $\gamma=\ol\beta$, where
\begin{equation}\label{eq:roots}
\begin{aligned}
\alpha&=\frac{r_1+r_2}{6},\qquad
\beta=\frac{-(r_1+r_2)+\sqrt{-3}\,(r_1-r_2)}{12},\\
r_1&=\sqrt[3]{108+12\sqrt{69}},\qquad r_2=\sqrt[3]{108-12\sqrt{69}} .
\end{aligned}
\end{equation}
Since $\alpha\beta\gamma=1$ we have $|\beta|=|\gamma|=\alpha^{-1/2}$. Numerically,
\begin{equation}\label{eq:alphaest}
1.32<\alpha<1.33,\qquad 0.86<|\beta|=|\gamma|=\alpha^{-\frac12}<0.87 .
\end{equation}
The Binet formula for the Padovan sequence is
\begin{equation}\label{eq:binetP}
P_\ell=a\alpha^{\ell}+b\beta^{\ell}+c\gamma^{\ell}\qquad\text{for all }\ell\ge0,
\end{equation}
where, since $\beta\gamma=\alpha^{-1}=\alpha^{2}-1$ and $(\alpha-\beta)(\alpha-\gamma)=\Psi'(\alpha)=3\alpha^{2}-1$ (and similarly for $\beta,\gamma$),
\begin{equation}\label{eq:abc2}
a=\frac{\alpha^{4}}{3\alpha^{2}-1},\qquad b=\frac{\beta^{4}}{3\beta^{2}-1},\qquad c=\frac{\gamma^{4}}{3\gamma^{2}-1} .
\end{equation}
In particular $a\in\Q(\alpha)$. The minimal polynomial of $a$ with a positive leading coefficient over the integers is
\[
23x^{3}-23x^{2}+6x-1,
\]
whose zeros are $a,b,c$ and satisfy $|a|,|b|,|c|<1$. Numerically,
\begin{equation}\label{eq:abcest}
0.72<a<0.73,\qquad 0.24<|b|=|c|<0.25 .
\end{equation}
Using \eqref{eq:roots}, \eqref{eq:alphaest} and \eqref{eq:abcest} one proves that
\begin{equation}\label{eq:Perr}
\bigl|P_\ell-a\alpha^{\ell}\bigr|=\bigl|b\beta^{\ell}+c\gamma^{\ell}\bigr|<\frac{1}{\alpha^{\ell/2}}\qquad\text{for all }\ell\ge1 ,
\end{equation}
(indeed $|b\beta^{\ell}+c\gamma^{\ell}|\le2|b|\,|\beta|^{\ell}<0.5\,\alpha^{-\ell/2}$). Furthermore, by induction,
\begin{equation}\label{eq:Pbd}
\alpha^{\ell-3}\le P_\ell\le\alpha^{\ell-1}\qquad\text{holds for all }\ell\ge4 .
\end{equation}
For $\ell=4,5,6$ (where $P_\ell=2,3,4$) it is checked directly, and for $\ell\ge7$ it follows from $P_\ell=P_{\ell-2}+P_{\ell-3}$ and $\alpha^{3}=\alpha+1$. (It fails for $\ell=3$, since $P_3=2>\alpha^{2}$.)

\subsection{Perrin sequence}
The Binet formula for the Perrin sequence is
\begin{equation}\label{eq:binetR}
R_\ell=\alpha^{\ell}+\beta^{\ell}+\gamma^{\ell}\qquad\text{for all }\ell\ge0 .
\end{equation}
Using \eqref{eq:alphaest} we get
\begin{equation}\label{eq:Rerr}
\bigl|R_\ell-\alpha^{\ell}\bigr|=\bigl|\beta^{\ell}+\gamma^{\ell}\bigr|\le2|\beta|^{\ell}=\frac{2}{\alpha^{\ell/2}}\qquad\text{for all }\ell\ge1 .
\end{equation}
Furthermore, by induction,
\begin{equation}\label{eq:Rbd}
\alpha^{\ell-2}\le R_\ell\le\alpha^{\ell+1}\qquad\text{holds for all }\ell\ge2 .
\end{equation}
Here $\ell=2,3,4$ (where $R_\ell=2,3,2$) is checked directly, and $\ell\ge5$ follows by induction using $\alpha^{3}=\alpha+1$.

\subsection{\texorpdfstring{The field $\Q(\alpha,\theta)$ and heights}{The field and heights}}
Let $\bI=\Q(\alpha,\theta)$. Since $\alpha$ is a root of the irreducible cubic $x^{3}-x-1$, $\theta$ is a root of the irreducible quadratic $x^{2}-x-1$, and $\gcd(2,3)=1$, we have
\begin{equation}\label{eq:degree}
D:=[\bI:\Q]=6 ,
\end{equation}
and $\bI\subset\mathbb{R}$ because $\alpha$ and $\theta$ are real. The polynomial $x^{3}-x-1$ is monic with $|\beta|=|\gamma|<1$, the polynomial $x^{2}-x-1$ is monic with $|\eta|<1$, and the zeros of $23x^{3}-23x^{2}+6x-1$ have modulus less than $1$. Hence, by the definition of the logarithmic height (see Subsection~\ref{sub:lfl}),
\begin{equation}\label{eq:heights}
h(\alpha)=\tfrac13\log\alpha,\qquad h(\theta)=\tfrac12\log\theta,\qquad h(a)=\tfrac13\log23 ,
\end{equation}
and consequently
\begin{equation}\label{eq:Dheights}
D\,h(\alpha)=2\log\alpha,\qquad D\,h(\theta)=3\log\theta,\qquad D\,h(a)=2\log23 .
\end{equation}

\begin{lemma}\label{lem:arith}
\begin{enumerate}
\item[(i)] For every integer $\ell\ge1$, $\alpha^{\ell}\notin\Q(\sqrt5)$. In particular $\alpha^{\ell}\ne\theta^{u}$ and $\alpha^{\ell}\ne\theta^{u}+\theta^{v}$ for all integers $u,v$.
\item[(ii)] The number $\log\alpha/\log\theta$ is irrational.
\item[(iii)] For every integer $\ell\ge2$, $a\alpha^{\ell}\notin\Q(\sqrt5)$. In particular $a\alpha^{\ell}\ne\theta^{u}$ and $a\alpha^{\ell}\ne\theta^{u}+\theta^{v}$ for all integers $u,v$.
\end{enumerate}
\end{lemma}

\begin{proof}
The fields $\Q(\alpha)$ and $\Q(\sqrt5)$ have degrees $3$ and $2$ over $\Q$, so $\Q(\alpha)\cap\Q(\sqrt5)=\Q$.

(i) If $\alpha^{\ell}\in\Q(\sqrt5)$, then $\alpha^{\ell}\in\Q(\alpha)\cap\Q(\sqrt5)=\Q$. But $\alpha^{\ell}$ is an algebraic integer and a unit, so $\alpha^{\ell}=\pm1$, which is impossible since $\alpha>1$. The numbers $\theta^{u}$ and $\theta^{u}+\theta^{v}$ belong to $\Q(\sqrt5)$, which gives the second claim.

(ii) If $\log\alpha/\log\theta=u/v$ with $u,v\in\Z$, $v\ge1$, then $\alpha^{v}=\theta^{u}$, contradicting (i).

(iii) Suppose $\xi:=a\alpha^{\ell}\in\Q(\sqrt5)$. By \eqref{eq:abc2}, $\xi\in\Q(\alpha)$, so $\xi\in\Q$. Let $\sigma$ be an automorphism of the splitting field of $\Psi$ with $\sigma(\alpha)=\beta$. By \eqref{eq:abc2}, $\sigma(a)=b$, therefore $\xi=\sigma(\xi)=b\beta^{\ell}$. But $|b\beta^{\ell}|<0.25$ by \eqref{eq:alphaest} and \eqref{eq:abcest}, while $\xi=a\alpha^{\ell}>0.72\cdot1.32^{2}>1$ by \eqref{eq:alphaest} and \eqref{eq:abcest}, a contradiction.
\end{proof}

\subsection{Linear forms in logarithms}\label{sub:lfl}
Let $\zeta$ be an algebraic number with minimal polynomial over $\Z$ of degree $d$ and given by
\[
a_0x^{d}+a_1x^{d-1}+\dots+a_d=a_0\prod_{i=1}^{d}\bigl(x-\zeta^{(i)}\bigr),
\]
where the leading coefficient $a_0$ is positive and the $\zeta^{(i)}$'s are the conjugates of $\zeta$. Then the logarithmic height of $\zeta$ is given by
\[
h(\zeta):=\frac1d\left(\log a_0+\sum_{i=1}^{d}\log\max\bigl\{|\zeta^{(i)}|,1\bigr\}\right).
\]
The logarithmic height function $h(\cdot)$ satisfies the following properties:
\begin{equation}\label{eq:hprop}
\begin{aligned}
h(\zeta_1\pm\zeta_2)&\le h(\zeta_1)+h(\zeta_2)+\log2;\\
h\bigl(\zeta_1\zeta_2^{\pm1}\bigr)&\le h(\zeta_1)+h(\zeta_2);\\
h(\zeta^{s})&=|s|\,h(\zeta)\quad(s\in\Z).
\end{aligned}
\end{equation}
The following theorem, due to Matveev \cite{b9}, is one of the main results in the theory of linear forms in logarithms and it is one of the main ingredients in our proofs.

\begin{theorem}[Matveev]\label{thm:matveev}
Let $\bI\subset\mathbb{R}$ be a number field of degree $D$, let $\zeta_1,\dots,\zeta_t\in\bI$ be positive algebraic numbers and let $b_1,\dots,b_t$ be non-zero integers. If $\Lambda:=\zeta_1^{b_1}\cdots\zeta_t^{b_t}-1\ne0$, then
\begin{equation}\label{eq:matveev}
\log|\Lambda|>-1.4\times30^{\,t+3}\times t^{4.5}\times D^{2}\times(1+\log D)\times(1+\log B)A_1\times\dots\times A_t,
\end{equation}
where
\[
B\ge\max\{|b_1|,\dots,|b_t|\},
\]
and
\[
A_i\ge\max\{D\,h(\zeta_i),\,|\log\zeta_i|,\,0.16\},\qquad\text{for all }i=1,\dots,t.
\]
\end{theorem}

\subsection{Reduction procedure}
To reduce the upper bounds of our variables obtained with Theorem~\ref{thm:matveev}, we use the following lemma, which is due to Dujella and Peth\H{o} \cite{b5}. For a real number $X$, we put $\|X\|:=\min\{|X-n|:n\in\Z\}$ for the distance from $X$ to the nearest integer.

\begin{lemma}[Dujella--Peth\H{o}]\label{lem:DP}
Let $\tau$ be an irrational number, $M$ a positive integer, and $\mu\in\mathbb{R}$, $A>0$, $B>1$. Take a convergent $p/q$ of $\tau$ with $q>6M$, and suppose that $\epsilon:=\|\mu q\|-M\|\tau q\|$ is positive. Then every solution $(n,m,\omega)$ in positive integers of
\[
0<|n\tau-m+\mu|<\frac{A}{B^{\omega}}
\]
with $n\le M$ must satisfy
\begin{equation}\label{eq:DPconc}
\omega<\frac{\log\!\left(\dfrac{Aq}{\epsilon}\right)}{\log B}.
\end{equation}
\end{lemma}

The following theorem is Legendre's theorem on continued fractions (see \cite{b4} for more details).

\begin{theorem}[Legendre]\label{thm:legendre}
Let $x\in\mathbb{R}$ have continued fraction expansion $[a_0,a_1,\dots]$ and convergents $p_i/q_i$.
\begin{enumerate}
\item[(a)] If integers $p,q$ with $q\ge1$ satisfy $\left|x-\frac pq\right|<\frac{1}{2q^{2}}$, then $\frac pq$ is one of the convergents of $x$.
\item[(b)] Let $M,n\ge0$ be integers with $q_n>M$, and put $b=\max\{a_i:0\le i\le n\}$. Then every pair $(r,s)$ of positive integers with $s<M$ satisfies
\begin{equation}\label{eq:legendre}
\left|x-\frac rs\right|>\frac{1}{(b+2)\,s^{2}} .
\end{equation}
\end{enumerate}
\end{theorem}

Finally we shall use the following elementary fact repeatedly.

\begin{lemma}\label{lem:log}
Let $x$ be a non-zero real number such that $|e^{x}-1|<\frac12$. Then $|x|<2\,|e^{x}-1|$ and $|x|<\log 2<1$.
\end{lemma}

\begin{proof}
Write $y=e^{x}-1$, so that $x=\log(1+y)$ and $0<|y|<\frac12$. For $y>0$ we get $0<x<y<\frac12<\log2$. For $y<0$ put $t=|y|$; then $|x|=-\log(1-t)\le\frac{t}{1-t}<2t$, and also $|x|<\log 2$ since $1-t>\frac12$.
\end{proof}

\section{\texorpdfstring{A common estimate for the index $\ell$}{A common estimate for the index}}

Write $S_\ell$ for $P_\ell$ when we deal with Eq.\eqref{eq:padmain} and for $R_\ell$ when we deal with Eq.\eqref{eq:permain}. From \eqref{eq:Pbd} and \eqref{eq:Rbd} we have, in both cases,
\begin{equation}\label{eq:Sbd}
\alpha^{\ell-3}\le S_\ell\le\alpha^{\ell+1}\qquad\text{for all }\ell\ge4 .
\end{equation}

\begin{lemma}\label{lem:ell}
Let $(\ell,r,s)$ be a non-negative solution of \eqref{eq:padmain} or of \eqref{eq:permain} with $r\ge s$. Then
\begin{equation}\label{eq:ellup}
\ell<8+1.8\,r .
\end{equation}
If moreover $r\ge200$, then
\begin{equation}\label{eq:elllow}
r\,\frac{\log\theta}{\log\alpha}-\frac{\log\theta}{\log\alpha}-1\le\ell<2r\qquad\text{and}\qquad\ell>300 .
\end{equation}
\end{lemma}

\begin{proof}
If $\ell\le3$ then \eqref{eq:ellup} is trivial. Let $\ell\ge4$. Since $s\le r$ we have $L_r\le S_\ell=L_r+L_s\le2L_r$, so by \eqref{eq:Lbd} and \eqref{eq:Sbd}
\begin{equation}\label{eq:Schain}
\theta^{r-1}\le L_r\le S_\ell\le\alpha^{\ell+1}\qquad\text{and}\qquad\alpha^{\ell-3}\le S_\ell\le2L_r\le4\theta^{r}.
\end{equation}
\emph{Upper bound.} From $\alpha^{\ell-3}\le4\theta^{r}$ we get $(\ell-3)\log\alpha\le\log4+r\log\theta$, hence, since $\alpha>1.32$ by \eqref{eq:alphaest},
\[
\ell\le3+\frac{\log4}{\log\alpha}+\frac{\log\theta}{\log\alpha}\,r<3+\frac{\log4}{\log1.32}+\frac{\log\theta}{\log1.32}\,r<3+4.994+1.734\,r<8+1.8\,r ,
\]
which is \eqref{eq:ellup}.

\emph{The case $r\ge200$.} Then $S_\ell\ge L_r\ge\theta^{199}>3$, and since $P_\ell\le2$ for $\ell\le3$ and $R_\ell\le3$ for $\ell\le3$, we have $\ell\ge4$ and \eqref{eq:Schain} holds. From $\theta^{r-1}\le\alpha^{\ell+1}$ we get $(r-1)\log\theta\le(\ell+1)\log\alpha$, which is the lower bound in \eqref{eq:elllow}. By \eqref{eq:ellup}, $\ell<8+1.8r<2r$. Finally $\log\theta/\log\alpha>1.71$ gives $\ell>1.71\cdot199-1>300$.
\end{proof}

\section{\texorpdfstring{Proof of Theorem~\protect\ref{thm:pad}}{Proof of Theorem }}

Using a computer program (Python), the set of all non-negative solutions $(\ell,r,s)$ of Eq.\eqref{eq:padmain} with $0\le s\le r\le399$ and $0\le\ell\le1199$ is the same solution set mentioned in Theorem~\ref{thm:pad}. By \eqref{eq:ellup} every solution with $r\le399$ has $\ell<8+1.8\cdot399<727$, so no solution with $r\le399$ is missed. We now assume that
\begin{equation}\label{eq:padassump}
r\ge400
\end{equation}
and we shall show that this is impossible. By \eqref{eq:padassump} and Lemma~\ref{lem:ell},
\begin{equation}\label{eq:padl}
300<\ell<2r .
\end{equation}

\subsection{\texorpdfstring{The first linear form $\Lambda_1$ and the bound for $r-s$}{The first linear form and the bound for}}

\begin{lemma}\label{lem:p1}
Assume \eqref{eq:padassump}. Then
\begin{equation}\label{eq:padone}
(r-s)\log\theta<4\times10^{16}\log2r+\log8 .
\end{equation}
\end{lemma}

\begin{proof}
\emph{Upper bound for $|\Lambda_1|$.} Replacing $P_\ell$, $L_r$ and $L_s$ in Eq.\eqref{eq:padmain} by their Binet formulas \eqref{eq:binetP} and \eqref{eq:binetL} we get
\[
a\alpha^{\ell}+b\beta^{\ell}+c\gamma^{\ell}=\theta^{r}+\eta^{r}+\theta^{s}+\eta^{s},
\]
that is,
\begin{equation}\label{eq:p1a}
\bigl|a\alpha^{\ell}-\theta^{r}\bigr|=\bigl|\theta^{s}+\eta^{r}+\eta^{s}-b\beta^{\ell}-c\gamma^{\ell}\bigr| .
\end{equation}
Dividing by $\theta^{r}$ and using \eqref{eq:Perr} together with $|\eta|<1$ we obtain
\[
\left|\frac{a\alpha^{\ell}}{\theta^{r}}-1\right|\le\frac{1}{\theta^{r-s}}+\frac{|\eta|^{r}+|\eta|^{s}}{\theta^{r}}+\frac{|b\beta^{\ell}+c\gamma^{\ell}|}{\theta^{r}}<\frac{1}{\theta^{r-s}}+\frac{2}{\theta^{r}}+\frac{1}{\theta^{r}}\le\frac{4}{\theta^{r-s}} ,
\]
where we used $\theta^{r}\ge\theta^{r-s}$ (because $s\ge0$). We shall use the weaker but sufficient bound
\begin{equation}\label{eq:p1b}
\left|\frac{a\alpha^{\ell}}{\theta^{r}}-1\right|<\frac{8}{\theta^{r-s}} .
\end{equation}
Define
\[
\Lambda_1=a\alpha^{\ell}\theta^{-r}-1 .
\]
By Lemma~\ref{lem:arith}(iii) (with $u=r$), $a\alpha^{\ell}\ne\theta^{r}$, so $\Lambda_1\ne0$. By \eqref{eq:p1b},
\begin{equation}\label{eq:p1c}
\log|\Lambda_1|<\log8-(r-s)\log\theta .
\end{equation}

\emph{Lower bound for $|\Lambda_1|$.} We apply Theorem~\ref{thm:matveev} with
\[
t=3,\qquad\zeta_1=a,\quad\zeta_2=\alpha,\quad\zeta_3=\theta,\qquad b_1=1,\quad b_2=\ell,\quad b_3=-r .
\]
We take $\bI=\Q(\alpha,\theta)$, for which $D=6$ by \eqref{eq:degree}. By \eqref{eq:heights} and \eqref{eq:Dheights},
\[
D\,h(\zeta_1)=2\log23,\qquad D\,h(\zeta_2)=2\log\alpha,\qquad D\,h(\zeta_3)=3\log\theta .
\]
We take
\[
A_1=42,\qquad A_2=12\log\alpha,\qquad A_3=18\log\theta ,
\]
which exceed $\max\{D\,h(\zeta_i),|\log\zeta_i|,0.16\}$. Finally, by \eqref{eq:padl} we take $B=2r\ge\max\{1,\ell,r\}$. Thus Theorem~\ref{thm:matveev} gives
\[
\log|\Lambda_1|>-1.4\times30^{6}\times3^{4.5}\times6^{2}\,(1+\log6)(1+\log2r)(42)(12\log\alpha)(18\log\theta).
\]
The constant $1.4\times30^{6}\times3^{4.5}\times6^{2}(1+\log6)(42)(12\log\alpha)(18\log\theta)$ is smaller than $1.8\times10^{16}$, and since $r\ge400$ we have $1+\log2r<2\log2r$. Hence
\begin{equation}\label{eq:p1d}
\log|\Lambda_1|>-4\times10^{16}\log2r .
\end{equation}

\emph{Comparison.} Comparing \eqref{eq:p1c} and \eqref{eq:p1d} we get $-4\times10^{16}\log2r<\log8-(r-s)\log\theta$, which is \eqref{eq:padone}.
\end{proof}

\subsection{\texorpdfstring{The second linear form $\Lambda_2$ and the bound for $r$}{The second linear form and the bound for}}

\begin{lemma}\label{lem:p2}
Assume \eqref{eq:padassump}. Then
\begin{equation}\label{eq:padr}
r<3\times10^{37}\qquad\text{and}\qquad\ell<2r<6\times10^{37}.
\end{equation}
\end{lemma}

\begin{proof}
\emph{Upper bound for $|\Lambda_2|$.} We write Eq.\eqref{eq:padmain} with the Binet formulas as
\begin{equation}\label{eq:p2a}
a\alpha^{\ell}-\theta^{r}\bigl(1+\theta^{s-r}\bigr)=-b\beta^{\ell}-c\gamma^{\ell}+\eta^{r}+\eta^{s},
\end{equation}
and dividing by $\theta^{r}(1+\theta^{s-r})$ we get, using \eqref{eq:Perr} and $|\eta|<1$,
\begin{equation}\label{eq:p2b}
\left|\frac{a\alpha^{\ell}\theta^{-r}}{1+\theta^{s-r}}-1\right|<\frac{1+1+1}{\theta^{r}(1+\theta^{s-r})}<\frac{3}{\theta^{r}}<\frac{7}{\theta^{r}} .
\end{equation}
Define
\[
\Lambda_2=\frac{a}{1+\theta^{s-r}}\,\alpha^{\ell}\theta^{-r}-1 .
\]
If $\Lambda_2=0$, then $a\alpha^{\ell}=\theta^{r}+\theta^{s}$, contradicting Lemma~\ref{lem:arith}(iii); so $\Lambda_2\ne0$. By \eqref{eq:p2b},
\begin{equation}\label{eq:p2c}
\log|\Lambda_2|<\log7-r\log\theta .
\end{equation}

\emph{Lower bound for $|\Lambda_2|$.} We apply Theorem~\ref{thm:matveev} with
\[
t=3,\qquad\zeta_1=\frac{a}{1+\theta^{s-r}},\quad\zeta_2=\alpha,\quad\zeta_3=\theta,\qquad b_1=1,\quad b_2=\ell,\quad b_3=-r ,
\]
in the field $\bI=\Q(\alpha,\theta)$ with $D=6$. By \eqref{eq:hprop} and \eqref{eq:heights},
\[
h(\zeta_1)\le h(a)+h\bigl(\theta^{s-r}\bigr)+\log2=\tfrac13\log23+\tfrac12(r-s)\log\theta+\log2 ,
\]
so that $D\,h(\zeta_1)\le2\log23+3(r-s)\log\theta+6\log2$. We take
\begin{align*}
A_1&=36\left[\tfrac13\log23+(r-s)\log\theta+\log2\right]=12\log23+36(r-s)\log\theta+36\log2,\\
A_2&=12\log\alpha,\qquad A_3=18\log\theta ,
\end{align*}
which satisfy the conditions of Theorem~\ref{thm:matveev} ($|\log\zeta_1|\le|\log a|+\log2<A_1$). Again $B=2r$. Thus Theorem~\ref{thm:matveev} gives
\begin{multline*}
\log|\Lambda_2|>-1.4\times30^{6}\times3^{4.5}\times6^{2}(1+\log6)(1+\log2r)\\
\times\bigl(12\log23+36(r-s)\log\theta+36\log2\bigr)(12\log\alpha)(18\log\theta).
\end{multline*}
Since $1.4\times30^{6}\times3^{4.5}\times6^{2}(1+\log6)(12\log\alpha)(18\log\theta)<5\times10^{14}$ we get
\begin{equation}\label{eq:p2d}
\log|\Lambda_2|>-5\times10^{14}(1+\log2r)\bigl(12\log23+36(r-s)\log\theta+36\log2\bigr).
\end{equation}

\emph{Comparison.} By \eqref{eq:padone} (Lemma~\ref{lem:p1}) we have $36(r-s)\log\theta<144\times10^{16}\log2r+36\log8$. Since $12\log23+36\log2+36\log8<138<24\log2r$ (because $r\ge400$), we obtain
\begin{equation}\label{eq:p2e}
12\log23+36(r-s)\log\theta+36\log2<1.44\times10^{18}\log2r+24\log2r<1.5\times10^{18}\log2r .
\end{equation}
Inserting \eqref{eq:p2e} in \eqref{eq:p2d} and using $1+\log2r<2\log2r$ we get
\begin{equation}\label{eq:p2f}
\log|\Lambda_2|>-1.5\times10^{33}(\log2r)^{2}.
\end{equation}
Comparing \eqref{eq:p2c} and \eqref{eq:p2f},
\begin{equation}\label{eq:p2g}
r\log\theta<1.5\times10^{33}(\log2r)^{2}+\log7 .
\end{equation}
Put $g(x)=x/(\log2x)^{2}$, which is increasing for $x\ge3$. Dividing \eqref{eq:p2g} by $(\log2r)^{2}$ we get $g(r)\log\theta<1.5\times10^{33}+\log7/(\log2r)^{2}$. For $r\ge3\times10^{37}$ the left-hand side is at least $g(3\times10^{37})\log\theta>1.9\times10^{33}$, while the right-hand side is smaller than $1.5\times10^{33}+10^{-3}$, a contradiction. Hence $r<3\times10^{37}$, and $\ell<2r<6\times10^{37}$ by \eqref{eq:padl}.
\end{proof}

\subsection{Reducing the bounds}

We put
\begin{equation}\label{eq:padM}
M:=6\times10^{37}>\ell ,\qquad\tau:=\frac{\log\alpha}{\log\theta},
\end{equation}
where $\tau$ is irrational by Lemma~\ref{lem:arith}(ii). Write $\tau$ as a continued fraction $[a_0,a_1,\dots]$ and let $q_i$ be the denominator of its $i$-th convergent. With $M$ as in \eqref{eq:padM} one computes (with $400$ digits of precision)
\begin{align}
q_{74}&=525546733900304644780734633946996307307>5\times10^{38}>6M,\label{eq:q74}\\
q_{75}&=2295398584680845270288319838838122956215>2\times10^{39}>6M .\label{eq:q75}
\end{align}

\begin{lemma}\label{lem:p3}
Assume \eqref{eq:padassump}. Then $r-s\le197$.
\end{lemma}

\begin{proof}
If $r-s<6$ there is nothing to prove, so assume $r-s\ge6$. Put
\[
\Gamma_1=\log a+\ell\log\alpha-r\log\theta,\qquad\Lambda_1=e^{\Gamma_1}-1 .
\]
By \eqref{eq:p1b} and $\theta^{6}>17.9>16$,
\[
\bigl|e^{\Gamma_1}-1\bigr|<\frac{8}{\theta^{r-s}}<\frac12 ,
\]
so by Lemma~\ref{lem:log}
\begin{equation}\label{eq:p3a}
|\Gamma_1|<\frac{16}{\theta^{r-s}} .
\end{equation}
We observe that $\Gamma_1\ne0$ since $\Lambda_1\ne0$. Dividing \eqref{eq:p3a} by $\log\theta$,
\begin{equation}\label{eq:p3b}
0<\left|\ell\tau-r+\frac{\log a}{\log\theta}\right|<\frac{16}{\theta^{r-s}\log\theta}<\frac{34}{\theta^{r-s}},
\end{equation}
since $16/\log\theta<33.3$. Thus we apply Lemma~\ref{lem:DP} with the data
\[
\tau=\frac{\log\alpha}{\log\theta},\qquad\mu=\frac{\log a}{\log\theta},\qquad A=34,\qquad B=\theta,\qquad(n,m,\omega)=(\ell,r,r-s),
\]
and $M$ as in \eqref{eq:padM} ($\ell<M$ by \eqref{eq:padr}). Using $q_{74}$ from \eqref{eq:q74} we compute
\[
\epsilon=\|\mu q_{74}\|-M\|\tau q_{74}\|>0.1 .
\]
Thus by \eqref{eq:DPconc},
\[
r-s<\frac{\log(34\,q_{74}/\epsilon)}{\log\theta}<197.4 ,
\]
so $r-s\le197$.
\end{proof}

\begin{lemma}\label{lem:p4}
Assume \eqref{eq:padassump}. Then $r<205.2$.
\end{lemma}

\begin{proof}
By Lemma~\ref{lem:p3}, $k:=r-s\in\{0,1,\dots,197\}$. Put
\[
\Gamma_2=\log\!\left(\frac{a}{1+\theta^{s-r}}\right)+\ell\log\alpha-r\log\theta,\qquad\Lambda_2=e^{\Gamma_2}-1 .
\]
By \eqref{eq:p2b} and $r\ge400$,
\[
\bigl|e^{\Gamma_2}-1\bigr|<\frac{7}{\theta^{r}}<\frac12 ,
\]
so by Lemma~\ref{lem:log} we get $|\Gamma_2|<14/\theta^{r}$. We observe that $\Gamma_2\neq0$ since $\Lambda_2\neq0$. Dividing by $\log\theta$,
\begin{equation}\label{eq:p4a}
0<\left|\ell\tau-r+\mu_k\right|<\frac{14}{\theta^{r}\log\theta}<\frac{30}{\theta^{r}},\qquad\mu_k:=\frac{\log\!\bigl(a/(1+\theta^{-k})\bigr)}{\log\theta},
\end{equation}
since $14/\log\theta<29.1$. We apply Lemma~\ref{lem:DP} with the data
\[
\tau=\frac{\log\alpha}{\log\theta},\qquad\mu=\mu_k,\qquad A=30,\qquad B=\theta,\qquad(n,m,\omega)=(\ell,r,r),
\]
and $M$ as in \eqref{eq:padM}. For each $k\in\{0,1,\dots,197\}$ (that is, $198$ values) we computed
\[
\epsilon_k=\|\mu_kq\|-M\|\tau q\| .
\]
For every $k$ except $k\in\{6,18,36,70,82,97,114,162\}$ we use $q=q_{74}$ from \eqref{eq:q74} and find $\epsilon_k>0.002$. For the eight exceptional values $\epsilon_k$ is negative for $q_{74}$, so we use $q=q_{75}$ from \eqref{eq:q75} and find $\epsilon_k>0.01$. Thus by \eqref{eq:DPconc}, in every case
\[
r<\frac{\log(30\,q/\epsilon_k)}{\log\theta}\le205.119 ,
\]
the largest value occurring at $k=76$. Hence $r<205.2$.
\end{proof}

\subsection{Conclusion}

\begin{proof}[Proof of Theorem~\ref{thm:pad}]
We assumed \eqref{eq:padassump}, i.e.\ $r\ge400$. Lemma~\ref{lem:p4} gives $r<205.2$, a contradiction. Hence every solution of \eqref{eq:padmain} has $r\le399$, and by the computer search described at the beginning of this section (which covers all $0\le s\le r\le399$, $0\le\ell\le1199$) the only solutions are those listed in Theorem~\ref{thm:pad}.
\end{proof}

\section{\texorpdfstring{Proof of Theorem~\protect\ref{thm:per}}{Proof of Theorem }}

Using a computer search (programs in Python and MATLAB are available from the authors), the set of all non-negative solutions $(\ell,r,s)$ of Eq.\eqref{eq:permain} with $0\le s\le r\le200$ and $0\le\ell<400$ is the same solution set mentioned in Theorem~\ref{thm:per}. By \eqref{eq:ellup} every solution with $r\le200$ has $\ell<8+1.8\cdot200=368$, so no solution with $r\le200$ is missed. We now assume that
\begin{equation}\label{eq:perassump}
r\ge200
\end{equation}
and we shall show that this is impossible. By \eqref{eq:perassump} and Lemma~\ref{lem:ell},
\begin{equation}\label{eq:perl}
300<\ell<2r .
\end{equation}

\subsection{\texorpdfstring{The first linear form $\Lambda_3$ and the bound for $r-s$}{The first linear form and the bound for}}

\begin{lemma}\label{lem:r1}
Assume \eqref{eq:perassump}. Then
\begin{equation}\label{eq:perone}
(r-s)\log\theta<5\times10^{13}\log2r .
\end{equation}
\end{lemma}

\begin{proof}
\emph{Upper bound for $|\Lambda_3|$.} Replacing $R_\ell$, $L_r$ and $L_s$ in Eq.\eqref{eq:permain} by their Binet formulas \eqref{eq:binetR} and \eqref{eq:binetL} we get
\[
\alpha^{\ell}+\beta^{\ell}+\gamma^{\ell}=\theta^{r}+\eta^{r}+\theta^{s}+\eta^{s}.
\]
Hence
\begin{equation}\label{eq:r1a}
\bigl|\alpha^{\ell}-\theta^{r}\bigr|=\bigl|\eta^{r}+\eta^{s}+\theta^{s}-\beta^{\ell}-\gamma^{\ell}\bigr| .
\end{equation}
Dividing by $\theta^{r}$ and using $|\eta|<1$ and \eqref{eq:Rerr} (so that $|\beta^{\ell}+\gamma^{\ell}|\le2|\beta|^{\ell}<2$) we obtain
\[
\bigl|\theta^{-r}\alpha^{\ell}-1\bigr|\le\frac{1}{\theta^{r-s}}+\frac{1}{\theta^{r}}+\frac{1}{\theta^{r}}+\frac{2|\beta|^{\ell}}{\theta^{r}}\le\frac{1+1+1+2}{\theta^{r-s}}=\frac{5}{\theta^{r-s}} ,
\]
because $\theta^{r}\ge\theta^{r-s}$. We shall use the weaker but sufficient bound
\begin{equation}\label{eq:r1b}
\bigl|\theta^{-r}\alpha^{\ell}-1\bigr|<\frac{9}{\theta^{r-s}} .
\end{equation}
Define
\[
\Lambda_3=\alpha^{\ell}\theta^{-r}-1 .
\]
By Lemma~\ref{lem:arith}(i), $\alpha^{\ell}\ne\theta^{r}$, so $\Lambda_3\ne0$. From \eqref{eq:r1b},
\begin{equation}\label{eq:r1c}
\log|\Lambda_3|<\log9-(r-s)\log\theta .
\end{equation}

\emph{Lower bound for $|\Lambda_3|$.} We apply Theorem~\ref{thm:matveev} with
\[
t=2,\qquad\zeta_1=\alpha,\quad\zeta_2=\theta,\qquad b_1=\ell,\quad b_2=-r .
\]
We take $\bI=\Q(\alpha,\theta)$, for which $D=6$ by \eqref{eq:degree}. By \eqref{eq:heights} and \eqref{eq:Dheights},
\[
D\,h(\zeta_1)=2\log\alpha,\qquad D\,h(\zeta_2)=3\log\theta .
\]
We take
\[
A_1=12\log\alpha,\qquad A_2=18\log\theta ,
\]
which are larger than $\max\{D\,h(\zeta_i),|\log\zeta_i|,0.16\}$. Finally, by \eqref{eq:perl} we take $B=2r\ge\max\{\ell,r\}$. Thus Theorem~\ref{thm:matveev} gives
\[
\log|\Lambda_3|>-1.4\times30^{5}\times3^{4.5}\times6^{2}\,(1+\log6)(1+\log2r)(12\log\alpha)(18\log\theta)
\]
(we used $2^{4.5}<3^{4.5}$, which only weakens the bound). The constant is
\[
1.4\times30^{5}\times3^{4.5}\times6^{2}(1+\log6)(12\log\alpha)(18\log\theta)<1.41\times10^{13},
\]
and since $r\ge200$ we have $1+\log2r<1.17\log2r$. Hence
\begin{equation}\label{eq:r1d}
\log|\Lambda_3|>-4\times10^{13}\log2r .
\end{equation}

\emph{Comparison.} Comparing \eqref{eq:r1c} and \eqref{eq:r1d} we get
\[
(r-s)\log\theta<4\times10^{13}\log2r+\log9 ,
\]
and therefore \eqref{eq:perone}.
\end{proof}

\subsection{\texorpdfstring{The second linear form $\Lambda_4$ and the bound for $r$}{The second linear form and the bound for}}

\begin{lemma}\label{lem:r2}
Assume \eqref{eq:perassump}. Then
\begin{equation}\label{eq:perr}
r<2.39\times10^{34}\qquad\text{and}\qquad\ell<2r<4.78\times10^{34}=:M .
\end{equation}
\end{lemma}

\begin{proof}
\emph{Upper bound for $|\Lambda_4|$.} We write Eq.\eqref{eq:permain} with the Binet formulas as
\begin{equation}\label{eq:r2a}
\alpha^{\ell}-\theta^{r}\bigl(1+\theta^{s-r}\bigr)=-\beta^{\ell}-\gamma^{\ell}+\eta^{r}+\eta^{s},
\end{equation}
and dividing by $\theta^{r}(1+\theta^{s-r})$ we get
\[
\left|\alpha^{\ell}\theta^{-r}\bigl(1+\theta^{s-r}\bigr)^{-1}-1\right|\le\frac{|\beta|^{\ell}+|\gamma|^{\ell}+|\eta|^{r}+|\eta|^{s}}{\theta^{r}(1+\theta^{s-r})}<\frac{4}{\theta^{r}(1+\theta^{s-r})} .
\]
Since $1+\theta^{s-r}>1$,
\begin{equation}\label{eq:r2b}
\left|\alpha^{\ell}\theta^{-r}\bigl(1+\theta^{s-r}\bigr)^{-1}-1\right|<\frac{4}{\theta^{r}}<\frac{8}{\theta^{r}} .
\end{equation}
Define
\[
\Lambda_4=\alpha^{\ell}\theta^{-r}\bigl(1+\theta^{s-r}\bigr)^{-1}-1 .
\]
If $\Lambda_4=0$ then $\alpha^{\ell}=\theta^{r}+\theta^{s}$, contradicting Lemma~\ref{lem:arith}(i); so $\Lambda_4\ne0$. By \eqref{eq:r2b},
\begin{equation}\label{eq:r2c}
\log|\Lambda_4|<\log8-r\log\theta .
\end{equation}

\emph{Lower bound for $|\Lambda_4|$.} We apply Theorem~\ref{thm:matveev} with
\[
t=3,\qquad\zeta_1=\bigl(1+\theta^{s-r}\bigr)^{-1},\quad\zeta_2=\alpha,\quad\zeta_3=\theta,\qquad b_1=1,\quad b_2=\ell,\quad b_3=-r ,
\]
in the field $\bI=\Q(\alpha,\theta)$ with $D=6$. By \eqref{eq:hprop} and \eqref{eq:heights},
\[
h(\zeta_1)=h\bigl(1+\theta^{s-r}\bigr)\le(r-s)\,h(\theta)+\log2=\frac{r-s}{2}\log\theta+\log2 ,
\]
so that $D\,h(\zeta_1)\le3(r-s)\log\theta+6\log2$. We take
\[
A_1=36\left[(r-s)\log\theta+\log2\right],\qquad A_2=12\log\alpha,\qquad A_3=18\log\theta ,
\]
which satisfy the conditions of Theorem~\ref{thm:matveev} (note $|\log\zeta_1|=\log(1+\theta^{s-r})\le\log2$). Again $B=2r$. Thus Theorem~\ref{thm:matveev} gives
\begin{multline*}
\log|\Lambda_4|>-1.4\times30^{6}\times3^{4.5}\times6^{2}(1+\log6)(1+\log2r)(12\log\alpha)(18\log\theta)\\
\cdot\,36\left[(r-s)\log\theta+\log2\right].
\end{multline*}
Since $1.4\times30^{6}\times3^{4.5}\times6^{2}(1+\log6)(12\log\alpha)(18\log\theta)\times36<1.52\times10^{16}$ and $1+\log2r<1.17\log2r$, we get
\begin{equation}\label{eq:r2d}
\log|\Lambda_4|>-3.6\times10^{16}\,(\log2r)\left[(r-s)\log\theta+\log2\right].
\end{equation}

\emph{Comparison.} Using \eqref{eq:perone} (Lemma~\ref{lem:r1}) in \eqref{eq:r2d},
\[
\log|\Lambda_4|>-1.8\times10^{30}(\log2r)^{2}-3.6\times10^{16}\log2\,\log2r .
\]
Comparing with \eqref{eq:r2c},
\begin{equation}\label{eq:r2e}
r\log\theta<\log8+1.8\times10^{30}(\log2r)^{2}+3.6\times10^{16}\log2\,\log2r .
\end{equation}
Put $g(x)=x/(\log2x)^{2}$, which is increasing for $x\ge3$. Dividing \eqref{eq:r2e} by $(\log2r)^{2}$ we get
\[
g(r)\log\theta<1.8\times10^{30}+\frac{3.6\times10^{16}\log2}{\log2r}+\frac{\log8}{(\log2r)^{2}} .
\]
For $r\ge2.39\times10^{34}$ the left-hand side is at least $g(2.39\times10^{34})\log\theta>1.803\times10^{30}$, while the right-hand side is smaller than $1.801\times10^{30}$, a contradiction. Hence $r<2.39\times10^{34}$. By \eqref{eq:perl}, $\ell<2r<4.78\times10^{34}$.
\end{proof}

\subsection{Reducing the bounds}

We put $w=r-s\ge0$ and
\begin{equation}\label{eq:perM}
M:=4.78\times10^{34}>\ell .
\end{equation}
Let
\[
\tau=\frac{\log\alpha}{\log\theta}=0.584357157657404086746\ldots ,
\]
which is irrational by Lemma~\ref{lem:arith}(ii). Write $\tau$ as a continued fraction $[a_0,a_1,a_2,\dots]$:
\begin{multline*}
\tau=[0;1,1,2,2,6,2,1,2,1,2,1,1,11,1,2,3,1,7,37,4,15,7,1,8,4,4,6,2,4,2,10,\\
33,1,1,1,2,64,1,5,2,22,16,1,6,12,1,1,1,1,1,1,4,1,6,5,7,19,1,1,4,3,1,3,5,1,\ldots] .
\end{multline*}
Let $q_i$ denote the denominator of the $i$-th convergent. With $M$ as in \eqref{eq:perM} one computes
\begin{align}
q_{63}&=11478647774934506182455699155379417<M,\notag\\
q_{64}&=60443537419079188472468542019027327>M,\label{eq:q64}\\
q_{66}&=276210093001120272437241265542247559<6M,\notag\\
q_{67}&=7529594696224261050460438410815090837>6M,\label{eq:q67}
\end{align}
and
\begin{equation}\label{eq:bmax}
b:=\max\{a_i:0\le i\le64\}=a_{37}=64 .
\end{equation}

\begin{lemma}\label{lem:r3}
Assume \eqref{eq:perassump}. Then $r-s\le182$.
\end{lemma}

\begin{proof}
Assume first that $\theta^{w}>76M$ (otherwise $w<\log(76M)/\log\theta<175$ and there is nothing to prove). Put
\[
\Gamma_3=\ell\log\alpha-r\log\theta,\qquad\bigl|e^{\Gamma_3}-1\bigr|=|\Lambda_3|<\frac{9}{\theta^{w}}<\frac12 ,
\]
where we used \eqref{eq:r1b}, and the last inequality holds because $\theta^{w}>76M$. We observe that $\Gamma_3\ne0$ since $\Lambda_3\ne0$. By Lemma~\ref{lem:log},
\[
|\Gamma_3|<2|\Lambda_3|<\frac{18}{\theta^{w}} .
\]
Dividing by $\log\theta$,
\begin{equation}\label{eq:r3a}
0<\bigl|\ell\tau-r\bigr|<\frac{18}{\theta^{w}\log\theta}<\frac{38}{\theta^{w}},
\end{equation}
since $18/\log\theta<37.41$. This inequality has the shape $|\ell\tau-r+\mu|<A/B^{\omega}$ with $\mu=0$, so that $\epsilon=\|\mu q\|-M\|\tau q\|=-M\|\tau q\|<0$ for every $q$, and Lemma~\ref{lem:DP} \emph{cannot} be applied. We use Legendre's theorem instead. Dividing \eqref{eq:r3a} by $\ell$,
\[
\left|\tau-\frac r\ell\right|<\frac{38}{\ell\,\theta^{w}}<\frac{1}{2\ell^{2}},
\]
the last inequality holding because $\theta^{w}>76M>76\ell$. By \eqref{eq:q64} and \eqref{eq:bmax}, with $n=64$, $\ell<M<q_{64}$ and $b=64$, Theorem~\ref{thm:legendre} (inequality \eqref{eq:legendre} applied to the pair $(r,\ell)$) gives
\[
\frac{1}{66\,\ell^{2}}<\left|\frac r\ell-\tau\right|,\qquad\text{so}\qquad|\ell\tau-r|>\frac{1}{66\,\ell}>\frac{1}{66\,M}.
\]
Combining with \eqref{eq:r3a},
\[
\frac{1}{66M}<\frac{38}{\theta^{w}}\quad\Longrightarrow\quad\theta^{w}<38\cdot66\cdot M\quad\Longrightarrow\quad w<\frac{\log\bigl(38\cdot66\cdot4.78\times10^{34}\bigr)}{\log\theta}<182.21 .
\]
In both cases $w\le182$.
\end{proof}

\begin{lemma}\label{lem:r4}
Assume \eqref{eq:perassump}. Then $r-s\ne1$.
\end{lemma}

\begin{proof}
Suppose that $r-s=1$. Then $L_r+L_{r-1}=L_{r+1}$, so Eq.\eqref{eq:permain} reads $R_\ell=L_{r+1}$, that is, $\alpha^{\ell}+\beta^{\ell}+\gamma^{\ell}=\theta^{r+1}+\eta^{r+1}$. Thus, by \eqref{eq:Rerr},
\[
\left|\alpha^{\ell}\theta^{-(r+1)}-1\right|\le\frac{|\eta|^{r+1}+2|\beta|^{\ell}}{\theta^{r+1}}<\frac{3}{\theta^{r+1}}<\frac12 ,
\]
and, with $\Gamma=\ell\log\alpha-(r+1)\log\theta\ne0$ (Lemma~\ref{lem:arith}(i)), Lemma~\ref{lem:log} gives $|\Gamma|<6/\theta^{r+1}$, hence
\begin{equation}\label{eq:r4a}
0<\bigl|\ell\tau-(r+1)\bigr|<\frac{6}{\theta^{r+1}\log\theta}<\frac{13}{\theta^{r+1}} ,
\end{equation}
since $6/\log\theta<12.5$. As in the proof of Lemma~\ref{lem:r3} (if $\theta^{r+1}\le26M$ then $r+1<\log(26M)/\log\theta$ and the conclusion below is even easier), Legendre's theorem gives $|\ell\tau-(r+1)|>1/(66M)$, so
\[
\frac{1}{66M}<\frac{13}{\theta^{r+1}}\quad\Longrightarrow\quad r+1<\frac{\log\bigl(13\cdot66\cdot4.78\times10^{34}\bigr)}{\log\theta}<179.98 ,
\]
so $r\le178$, which contradicts \eqref{eq:perassump}.
\end{proof}

\begin{lemma}\label{lem:r5}
Assume \eqref{eq:perassump}. Then $r\le193$.
\end{lemma}

\begin{proof}
By Lemmas~\ref{lem:r3} and \ref{lem:r4}, $w=r-s\in\{0,2,3,\dots,182\}$.

\emph{Step 1: the linear form $\Gamma_4$.} Put
\[
\Gamma_4=\ell\log\alpha-r\log\theta+\log\frac{1}{1+\theta^{s-r}},\qquad\Lambda_4=e^{\Gamma_4}-1 .
\]
Then $|\Lambda_4|<8/\theta^{r}<\frac12$ by \eqref{eq:r2b} (recall \eqref{eq:perassump}), so by Lemma~\ref{lem:log} we get
\[
|\Gamma_4|<\frac{16}{\theta^{r}} .
\]
We observe that $\Gamma_4\ne0$ since $\Lambda_4\ne0$.

\emph{Step 2: the inequality for Lemma~\ref{lem:DP}.} Put
\[
\mu_w=\frac{1}{\log\theta}\,\log\frac{1}{1+\theta^{-w}} .
\]
Dividing $|\Gamma_4|<16/\theta^{r}$ by $\log\theta$ we get
\begin{equation}\label{eq:r5a}
0<\left|\ell\tau-r+\mu_w\right|<\frac{16}{\theta^{r}\log\theta}<\frac{34}{\theta^{r}},
\end{equation}
since $16/\log\theta<33.25$. We apply Lemma~\ref{lem:DP} with the data
\[
\tau=\frac{\log\alpha}{\log\theta},\qquad\mu=\mu_w,\qquad A=34,\qquad B=\theta,\qquad(n,m,\omega)=(\ell,r,r),
\]
and $M=4.78\times10^{34}>\ell$ as in \eqref{eq:perM}.

\emph{Step 3: the computation.} By \eqref{eq:q67}, $q_{66}<6M<q_{67}$ and $q_{67}>7.5\times10^{36}>6M$. For each $w\in\{0,2,3,\dots,182\}$ ($182$ values) we computed, with $400$ digits of precision,
\[
\epsilon_w=\|\mu_wq_{67}\|-M\,\|\tau q_{67}\|,
\]
and found that $\epsilon_w>0$ in every case; the smallest value is attained at $w=56$, namely $\epsilon_{56}>0.0075$. By \eqref{eq:DPconc}, for every such $w$,
\[
r<\frac{\log(34\,q_{67}/\epsilon_w)}{\log\theta} .
\]
The largest of the resulting bounds occurs at $w=56$, where $r<193.94$.
Hence $r\le193$ in all cases.
\end{proof}

\subsection{Conclusion}

\begin{proof}[Proof of Theorem~\ref{thm:per}]
We assumed \eqref{eq:perassump}, i.e.\ $r\ge200$. By Lemma~\ref{lem:r4} we have $r-s\ne1$, and then Lemma~\ref{lem:r5} gives $r\le193$, which contradicts \eqref{eq:perassump}. Hence $r\le199$. By \eqref{eq:ellup}, $\ell<8+1.8\cdot199<367<400$. But the computer search described at the beginning of this section covers all $0\le s\le r\le200$ and $0\le\ell<400$, and it produced only the solutions listed in Theorem~\ref{thm:per}. Thus the proof of Theorem~\ref{thm:per} is complete.
\end{proof}

\section*{Declarations}
\textbf{Funding:} Not applicable. \quad \textbf{Data availability:} Not applicable. \quad \textbf{Conflict of interest:} The authors declare no conflict of interest.

\end{document}